\documentstyle[12pt,amssymb]{amsart}

\newcommand{\cl}{\operatorname{cl}}
\newcommand{\inte}{\operatorname{int}}

\newtheorem{theor}{Theorem}
\newtheorem{lemma}{Lemma}

\newtheorem{defin}{Definition}
\newtheorem{fact}{Fact}

\newtheorem{claim}{Claim}

\renewcommand{\thehip}

\theoremstyle{definition}

\newtheorem{rem}{Remark}

\newtheorem{rems}{Remarks}

\theoremstyle{remark}

\begin{document}

\begin{center}
{\bf REGULAR SUBALGEBRAS OF COMPLETE  \\ BOOLEAN ALGEBRAS}
\end{center}
\bigskip

\begin{center}
{\sc Aleksander B{\l}aszczyk and Saharon Shelah}~\footnote{The
research of the second author was partially supported by the
Basic Research Foundation of the Israel Academy of Sciences and
Humanities. This publication has Number 640 in S. Shelah's
list. \newline November 17, 1997}
\end{center} 


\null\hfill\parbox[t]{11cm}{{\small\bf Abstract.} \small There
is shown that there exists a complete, atomless,
$\sigma$-centered Boolean algebra, which does not contain any
regular countable subalgebra if and only if there exist a
nowhere dense ultrafilter. Therefore the existence of such
algebras is undecidable in ZFC.} \hfill \null 
\bigskip
\vspace{1 cm}

A subalgebra $\Bbb{B}$ of a Boolean algebra $\Bbb{A}$ is called
regular whenever for every $X\subseteq \Bbb{B}$,
$\sup_{\Bbb{B}}X=\bold{1}$ implies $\sup_{\Bbb{A}}X=\bold{1}$;
see e.g. Heindorf and Shapiro \cite{HS}.  Clearly, every dense
subalgebra is regular. Although, every complete Boolean algebra
contains a free Boolean algebra of the same size (see the
Balcar-Franek Theorem; \cite{BF}), not always such an embedding
is regular.  For instance, if $\Bbb{B}$ is a measure algebra,
then it contains a free subalgebra of the same cardinality as
$\Bbb{B}$, but $\Bbb{B}$ cannot contain any free Boolean algebra
as a regular subalgebra. Indeed, measure algebras are weakly
$\sigma$-distributive and free Boolean algebra are not and a
regular subalgebra of a $\sigma$-distributive one is again
$\sigma$-distributive. Thus $\Bbb{B}$ does not contain any free
Boolean algebra.  On the other hand, measure algebras are not
$\sigma$-centered. So, a natural question arises whether there
can exists a $\sigma$-centered, complete, atomless Boolean
algebra $\Bbb{B}$ without regular free subalgebras. Since
countable atomless Boolean algebras are free and every free
Boolean algebra contains a countable regular free subalgebra, it
is enough to ask whether $\Bbb{B}$ contains a countable regular
subalgebra. In the paper we prove that such an algebra exists
iff there exists a nowhere dense ultrafilter.

\begin{defin}[Baumgartner Principle, see \cite{Bau}] A filter
$D$ on $\omega$ is called nowhere dense if for every function
$f$ from $\omega$ to the Cantor set ${}^\omega 2$ there exists a
set $A\in D$ such that $f(A)$ is nowhere dense in ${}^\omega 2$.
\end{defin} 

In the sequel we will rather interested in nowhere dense
ultrafilters. Observe that every $P$-ultrafilter (i.e. every
$P$-point in $\omega^{*}$) is a nowhere dense ultrafilter.

\begin{theor}
There exists an atomless, complete, $\sigma$-centered Boolean algebra
without any countable regular subalgebras iff there exists a
nowhere dense ultrafilter.
\end{theor}

In the first part of the paper there are used forcing methods to
show that a nowhere dense ultrafilters exist whenever there exists a
$\sigma$-centered forcing $\Bbb{P}$ such that above every
element of $\Bbb{P}$ there are two incompatible ones and
$\Bbb{P}$ does not add any Cohen real. The forcing constructed
here uses some ideas from Gitik and Shelah \cite{GS}. They have
shown that if $\Bbb{P}$ is a 
$\sigma$-centered forcing notions and $\{A_n\colon n<\omega\}$
are subsets of $\Bbb{P}$ witnessing this and both $\Bbb{P}$ and
$A_n$'s are Borel, then $\Bbb{P}$ adds a Cohen real.  On the
other hand there is known that a forcing
$\Bbb{P}$ adds a Cohen real iff the complete Boolean algebra
$ \Bbb B = RO(\Bbb P) $ being the completion of $\Bbb{P}$
contains an element $ u $ such that the rediduced Boolean
algebra $ \Bbb B | u $ has a regular free
Boolean subalgebra. Thus, to prove the Theorem 1 we need to show
in particular the following:     

\begin{theor}
If there exists a $\sigma$-centered forcing $\Bbb{P}$ such that above
every element of $\Bbb{P}$ there are two incompatible ones and
$\Bbb{P}$ does not add any Cohen real then there exists a nowhere dense
ultrafilter on $\omega$.
\end{theor}

We shall precede the proof by some definitions and a lemma.

\begin{defin}
{\rm (a)} \ A forcing $\Bbb{P}$ is called $\sigma$-centered if
${\Bbb{P}}=\bigcup\{A_n\colon n<\omega\}$ where each $A_n$ is directed i.~e.
for every $p,q\in A_n$ there exists $r\in A_n$ such that $p\leqslant r$ and
$q\leqslant r$.

{\rm (b)} \ A forcing $ \Bbb P $ adds a Cohen real if there exists a $
\Bbb P $--name $ \underline{\underline{r}} \in ^{\omega} 2 $
such that for every open dense set $ {\cal D} \subset ^{\omega}
2 $ there is $ \Vdash_{\Bbb P} `` \underline{\underline{r}} \in
{\cal D}^*" $, where $ {\cal D}^* $ denotes the encoding if $
{\cal D} $ in the Boolean universe. 
\end{defin}

\begin{rems}\

(a)
The order of forcing in this notation is inverse of the one in
the Boolean algebra.

(b)
We can just assume that there is a member $p$ of
$\Bbb{P}$ such that if $q$ is above $p$ then 
there are $r_1$ and $r_2$ above
$q$ which are incompatible in $\Bbb{P}$.

\end{rems}

\begin{defin}
A set $X\subseteq {}^{\omega>}2$ is somewhere dense if there exists an
$\eta\in {}^{\omega>}2$ such that for every $\nu\in{}^{\omega>}2$ there
is  $\varrho\in X$ with $\eta\sphat\: \nu\unlhd \varrho$, where
$\eta\sphat\:\nu$ stands for the concatenation of $\eta$ and $\nu$ and the
relation $\unlhd$ means that $\varrho$ is an extension of the sequence
$\eta\sphat\:\nu$.
\end{defin}

\begin{lemma}
A filter $D$ on $\omega$ is not nowhere dense iff it is a
so-called well behaved filter, i.e. there is a function $f\colon
\omega\to{}^{\omega>}2$ such that for every $B\in D$ the range
of $f$ restricted to $B$ is somewhere dense.
\end{lemma}

\begin{pf}
Suppose $f\colon \omega\to{}^{\omega>}2$ and $B\subseteq \omega$
and the image of $B$ is not nowhere dense. Without loss of
generality we can assume that the range of $f$ is dense in
itself. Since every closed and dense in itself subset of the
Cantor cube ${}^{\omega}2$ is homeomorphic to the whole
${}^{\omega}2$ we can assume also that the range of $f$ is dense
in ${}^{\omega}2$. Moreover, since it is countable it can be
identified with a subset of the set ${}^{\omega>}2$ of all
rational points of the Cantor set. Thus without loss of
generality we can assume that $f$ maps $\omega$ into
${}^{\omega>}2$. On the other hand a set
$X\subseteq{}^{\omega>}2$ is nowhere dense whenever for every
$\eta\in {}^{\omega>}2$ there exists some $\nu\in{}^{\omega>}2$
such that the set of all sequences extending $\eta\sphat\:\nu$
is disjoint with $X$. Therefore, since the image of $B$ under
$f$ is not nowhere dense in ${}^{\omega>}2$, it can be
identified with a somewhere dense subset of ${}^{\omega>}2$.
This in fact completes the proof of the lemma.
\end{pf}

\begin{rem}
If $D$ is a filter on $\omega$ and $\cal{P}(\omega)/D$ is
infinite then $D$ is not nowhere dense. Indeed, if $\langle
A_n\colon n<\omega\rangle$ is a partition of $\omega$ such that
$\omega\setminus A_n\notin D$ for all $n<\omega$ and $\langle
e_n\colon n<\omega\rangle$ list the set ${}^{\omega>}2$ then the
map $f\colon \omega\to{}^{\omega>}2$ defined by the formula
$$
f(e)=e_n\quad\text{iff}\quad e\in A_n
$$ 
witness ``$D$ is well behaved".
\end{rem}

\begin{pf}[of Theorem 2
]
Assume that there are no nowhere dense ultrafilter. Further
assume that $\Bbb{P}$ is a forcing in which above each element
there are two incompatible one and ${\Bbb{P}}=\bigcup\{A_n\colon
n<\omega\}$ where each $A_n$ is directed. We start with the
following known fact which we prove here for the sake of
completeness: 

\begin{fact}[0]
Forcing with ${\Bbb{P}}$ add a new real.
\end{fact}

In fact, by assumption forcing with $\Bbb{P}$
add a new subset to $\Bbb{P}$, hence a new subset to some
ordinal. In the set 
\begin{multline*}
{\cal{K}}=\{(\alpha,p,\underset{\sptilde{}}{\tau})\colon
p\in{\Bbb{P}},\alpha\;{\text{an ordinal
and}}\;\underset{\sptilde{}}{\tau}\\{\text{a}}
\;{\Bbb{P}}-{\text{name of a subset of
}}\alpha{\text{ such that
}}p\Vdash``\underset{\sptilde{}}{\tau}\notin V"\} 
\end{multline*}
we choose $(\alpha,p,\underset{\sptilde{}}{\tau})$ with $\alpha$
being minimal. So necessarily $\alpha$ is  a cardinal and

$$p\Vdash ``\text{the tree } ({}^{\alpha>}2,\unlhd)\;\text{has a new
}\alpha\text{-branch in }V^{\Bbb{P}}"$$

So, as $\Bbb{P}$ satisfies ccc, necessarily
$\text{cf}(\alpha)=\aleph_0$ and letting
$\alpha=\bigcup_{n<\omega}\alpha_n$, where
$\alpha_n<\alpha_{n+1}$ 

for some countable $w\subseteq {}^{\alpha>}2$ we get
$$
p\Vdash``(\forall
n<\omega)(\underset{\sptilde{}}{\tau}\restriction \alpha_n\in
w)", $$
so $p\Vdash``\text{we add a new subset to }w, |w|=\aleph_0"$, as required.

Now we can restrict ourselves to $\Bbb{P}$ as above or repeat the argument
above any $q\in\Bbb{P}$.

Now we fix a  $\Bbb{P}$-name of a new real
$\underline{\underline{r}}\in{}^{\omega}2$ added by $\Bbb{P}$.
For every $p\in \Bbb{P}$ we set
$$
T_p=\{\eta\in{}^{\omega>}2\colon \neg (p\Vdash\neg(``\eta\trianglelefteq
\underline{\underline{r}}"))\}
$$
i.e. $\eta\in T_p$ iff there exists $q\in\Bbb{P}$ such that $p\leqslant q$
and $q\Vdash``\eta=\underline{\underline{r}}\restriction\;
\text{lg}\:\eta"$,
where $\text{lg}\:\eta$ denotes the length of the sequence
$\eta$.
\begin{fact}[1]
For every $p\in\Bbb{P}$, $T_p$ is a subtree of ${}^{\omega>}2$, i.e
$\eta\unlhd \nu$ and $\nu\in T_p$ implies $\eta\in T_p$ and
$\langle\rangle\in T_p$ , where $\langle\rangle$ denotes the empty
sequence.
\end{fact}
Indeed, if $\eta\unlhd \nu$ and
$\nu=\underline{\underline{r}}\restriction\; \text{lg}\:\nu$, then
$\eta=\underline{\underline{r}}\restriction\; \text{lg}\:\eta$.

\begin{fact}[2]
The tree $T_p$ has no maximal elements.
\end{fact}
To prove the fact (2) we fix $\eta\in T_p$. Then there is $q\in\Bbb{P}$
such that $p\leqslant q$ and
$$q\Vdash ``\underline{\underline{r}}\restriction\;\text{lg}\:\eta=\eta".$$
Let $k=\text{lg}(\eta)$, so $I=\{r\in{\Bbb{P}}\colon r\;\text{forces a
value to }\underline{\underline{r}}\restriction(k+1)\}$ is a dense and
open subset of $\Bbb{P}$, hence there is $q'\in{\Bbb{P}}$ such that $q\leqslant q'$ and
$q'$ forces a value to $\underline{\underline{r}}\restriction(k+1)$, say $\vartheta$. So $q'$
also forces $\underline{\underline{r}}\restriction k=\vartheta\restriction k$, but $q\leqslant q'$ and
$q\Vdash ``\underline{\underline{r}}\restriction k=\eta\;\text{hence
}\vartheta\restriction k=\eta"$. As $q'$ witness $\vartheta\in T_p$ and $\vartheta\in {}^{k+1}2$
and $\eta\in {}^k2$, $\eta\unlhd \vartheta$ which completes the proof of fact (2).

\begin{fact}[3]
The set $\lim T_p$ of all $\omega$-branches is closed, i.e. if
 $\eta\in{}^ {\omega}2\setminus\lim T_p$ then there exists
 $\nu\in{}^{\omega>}2$ such that $\nu\unlhd \eta$ and the set of all
 $\omega$-branches extended $\eta$ is disjoint with $\lim T_p$.
\end{fact}
Indeed, if $\eta\in {}^\omega2\setminus \lim T_p$ then there exists $n\in\omega$ such that
$n\leqslant m<\omega$ implies $\eta\restriction m\notin T_p$. 
By Fact 1 it is clear that every $\omega$-branch extending
$\nu$ does not belong to $T_p$, which proves the Fact 3.

Now let us observe that the family
$$\{T_p\colon p\in A_n\}$$
is directed under inclusion, i.e. if $p,q\in A_n$ and $r\in\Bbb{P}$ is
such that $p\leqslant r$ and $q\leqslant r$ then
$$T_r\subseteq T_p\cap T_q.$$
Indeed, if $\eta\in{}^{\omega>}2$ and there exists $s\geqslant r$ such
that $s\Vdash``\eta=\underline{\underline{r}}\restriction\text{lg}\:\eta"$
then of course $s\geqslant p$ and $s\geqslant q$ and thus $\eta$ belongs
to $T_p$ and $T_q$.

So by compactness of ${}^{\omega}2$ and Facts 1-3 we get the following:

\begin{fact}[4]
The set
$$T_n=\bigcap\{T_p\colon p\in A_n\}$$
is a subtree of ${}^{\omega>}2$ and the set of $\omega$-branches of $T_n$
is non-empty.  \end{fact}

Now we make a choice:
\begin{equation}
\eta_n^{*}\quad\text{is an}\;\omega\;\text{- branch of}\;T_n.
\end{equation}

Subsequently for every $n<\omega$ and every $p\in A_n$ we define
$$B_p^n=\{k<\omega\colon(\exists q)(q\in{\Bbb{P}})(p\leqslant q\land
q\Vdash``\underline{\underline{r}}\restriction k=\eta^*_n\restriction
k\And\underline{\underline{r}}(k)\ne\eta^*_n(k)")\}$$
We have the following:
\begin{fact}[5]
For every $n<\omega$ and every $p\in A_n$ the set $B_p^n$ is infinite.
\end{fact}

Indeed, since $p\in A_n$ and $T_n$ is a subtree of $T_p$, $\eta^*_n$ is an
$\omega$-branch of $T_p$. Let us fix $m<\omega$. Then, by the definition
of $T_p$, there exists $r\in\Bbb{P}$ such that $r\geqslant p$ and
$$r\Vdash``\eta^*_n\restriction m=\underline{\underline{r}}\restriction
m".$$
On the other hand
$$\Vdash_{\Bbb{P}}``\underline{\underline{r}}\ne\eta_n^* ",$$
because $\underline{\underline{r}}$ is a new real. Thus for some
$q\in\Bbb{P}$, $q\geqslant r$ and $k<\omega$ we get
$$q\Vdash ``\underline{\underline{r}}\restriction k\ne \eta_n^*\restriction
k".$$
We can assume that $k$ is minimal with such a property. Since $r\leqslant
q$, it must be $k>m$. But $q\geqslant p$ and thus, by minimality of $k$,
we have $k-1\in B_p^n$, which proves the Fact 5.

Now we establish for every $n<\omega$ the following definition:
$${\cal{D}}_n^0=\{B\subseteq\omega\colon(\exists p)(p\in
A_n)(|B_p^n\setminus B|<\omega)\}.$$

\begin{fact}[6]
For every $n<\omega$, ${\cal{D}}_n^0$ is a filter. 
\end{fact}

Indeed, there exists $p_1,p_2\in A_n$ such that both $B_{p_1}^n\setminus
 B_1$ and $B_{p_2}^n\setminus B_2$ are finite. Since $A_n$ is directed we
 can choose $r\in A_n$ such that $p_1\leqslant r$ and $p_2\leqslant r$. On
 the other hand the definition of $B_p^n$ easily follows that
 $$ p\leqslant q\quad\text{implies}\quad B_q^n \subseteq B_p^n.$$
Thus $B_r^n\subseteq B_{p_1}^n\cap B_{p_2}^n$ and therefore
$$B_r^n\setminus(B_1\cap B_2)\subseteq (B_{p_1}^n\setminus
 B_1)\cup(B_{p_2}^n\setminus B_2)$$ is finite. Clearly, 
every superset of an element of ${\cal{D}}_n^0$ also belongs to
${\cal{D}}_n^0$ and,
by the Fact 5,
${\cal{D}}_n^0$ does not contain empty set, which completes the proof of
Fact 6.

Now by Fact 5 and Fact 6, we can make the following choice: for $n<\omega$
\begin{equation}
{\cal{D}}_n\;\text{is a non-principial ultrafilter
                containing}\;{\cal{D}}_n^0
\end{equation}

By our hypothesis the ultrafilters ${\cal{D}}_n$ are not nowhere dense and so
by Lemma (3.~1.~D) for every $n<\omega$ we can choose a function $f_n\colon
\omega\to {}^{\omega>}2$ such that
\begin{equation}
(\forall B\in {\cal{D}}_n)(\exists u\in{}^{\omega>}2)(\forall \nu\in
{}^{\omega>}2)(\exists k\in B)(u\sphat\:\nu\unlhd f_n(k)).
\end{equation}

Without loss of generality we may assume that the empty sequence does not
belong to the range of $f_n$.

Now we have to come back to the sequence $\{\eta_n^*\colon n<\omega\}$ of
$\omega$-branches of the trees $T_n$. Since it can happen that the
sequence is not one-to-one we consider the set
$$Y=\{n<\omega\colon\eta_n^*\notin\{\eta_m^*\colon m<n\}\}.$$
Then for $n,m\in Y$ we have $\eta_n^*\ne\eta_m^*$ whenever $n\ne m$.

In the sequel we shall need the following:

\begin{claim}
If $\langle \eta_n\colon n<\omega\rangle\subseteq {}^{\omega}2$ is a
        sequence of pairwise different $\omega$-branches of a tree $T\subseteq
        {}^{\omega>}2$ there exists an increasing sequence $\langle e_n\colon
        n<\omega\rangle\subseteq\omega$ such that for all $n<m<\omega$ we have
\begin{equation}
\{\eta_n\restriction l\colon e_n<l<\omega\}\cap\{\eta_m\restriction
        l\colon e_m<l<\omega\}=\emptyset.
\tag{$*$}\end{equation}
\end{claim}
To prove the claim observe that $\eta_n\restriction l\ne\eta_m\restriction
        l$ and $k>l$ implies $\eta_n\restriction k\ne\eta_m\restriction
        k$. Now assume that $e_0,\dots,e_n$ are defined so that the
        condition ($*$) holds true. Since
        $\eta_{n+1}\notin\{\eta_0,\dots,\eta_n\}$ there exists $k<\omega$
        such that $\eta_0\restriction k,\dots,\eta_n\restriction
        k,\eta_{n+1}\restriction k$ are pairwise different. We can assume that
        $k>e_n$ and $e_{n+1}$ to be the first such $k$. This completes the
proof of the claim.

Now using the Claim we can choose an increasing sequence $\langle
        e_n\colon n<\omega\rangle\subseteq \omega$ in such a way that letting
 $$C_n=\{\eta_n^*\restriction l\colon e_n\leqslant l<\omega\}.$$
the sequence $\langle C_n\colon n\in Y\rangle$ consists of pairwise
        disjoint sets and we have
        $$\eta_n^*=\eta_m^*\Leftrightarrow e_n=e_m\Leftrightarrow C_n=C_m.$$

        For each $(n,m)\in\omega\times\omega$ let us choose a set
        $A_{n,m}\in{\cal{D}}_n$ in such a way that
        $${\cal{D}}_n\ne{\cal{D}}_m\Rightarrow \omega\setminus
        A_{n,m}\in{\cal{D}}_m.$$
We also define a $\Bbb{P}$-name $\underset{\sptilde{}}{\tau}@!_{n,m}$ as
follows: for $G\subseteq\Bbb{P}$ generic over $V$ we set
$$\underset{\sptilde{}}{\tau}@!_{n,m}[G]=
\begin{cases}
1 & \text{if }\text{lg}@,(\eta^*_n\cap \underline{\underline{r}})
\in A_{n,m},\\ 0 & \text{otherwise} \end{cases}$$
where for $\eta,\nu\in{}^{\omega}2$ the symbol $\eta\cap\nu$ denotes
the longest common initial segment of $\eta$ and $\nu$.

Now let $\langle(m_i,n_i)\colon i<\omega\rangle$ list $\omega\times\omega$
and let define a $\Bbb{P}$-name $\underline{\underline{r}}^\otimes$ of
a member of ${}^{\omega}2$ as follows:
$$\underline{\underline{r}}^\otimes(2i+1)=
\underset{\sptilde{}}{\tau}@!_{n_i,m_i}\quad\text{and}\quad
\underline{\underline{r}}^\otimes(2i)=\underline{\underline{r}}(i).$$
Then for every $p\in\Bbb{P}$ we define
$$T_p^{\otimes}=\{\nu\in{}^{\omega>}2\colon \neg(p\Vdash\neg(\nu\unlhd
\underline{\underline{r}}^\otimes))\}.$$

Similarly like in Facts 1-3 one can show that
for every $p\in\Bbb{P}$, $T_p^{\otimes}$
is a subtree of ${}^{\omega>}2$  with no maximal elements with the set
$\lim T_p^{\otimes}$ of all $\omega$-branches being closed. Then,
similarly like in the Fact~4, for $n<\omega$,
$$T_n^\otimes=\bigcap\{T_p^\otimes\colon p\in A_n\}$$ is a non-empty tree
contained in ${}^{\omega>}2$. Those trees have the following property:

\begin{equation}
\text{for }n<\omega\;\text{and }p\in A_n\;\text{we have  }
T_p=\{\langle\nu(2i)\colon
2i<\text{lg}\,\nu\rangle\colon\nu\in T_p^\otimes\},
\tag{$*_1$}\end{equation}
As a consequence we get:
\begin{equation}
\text{for } n<\omega\;\text{we have }T_n=\{\langle\nu(2i)\colon
2i<\text{lg}\,\nu\rangle\colon\nu\in T_p^n\}.
\tag{$*_2$}\end{equation}

\begin{fact}[7]
For every $n<\omega$ we can find an $\omega$-branch
$\eta^\otimes_n\in\lim T_n^\otimes$ such that:

a) $\eta^*_n(i)=\eta_n^\otimes(2i)$ for $i<\omega$.

b) If $i<\omega$, then $\eta^{\otimes}_n(2i+1)$ is 1 if $A_{m_i,n_i}\in{\cal{D}}_n$ and it equals 0 if
$A_{m_i,n_i}\notin{\cal{D}}_n$.
\end{fact}

Hence by condition ($*_1$) we get
\begin{equation}
\text{if }m<n<\omega\;\text{and }\eta_n^*=\eta_m^*\;\text{and
}{\cal{D}}_n\ne{\cal{D}}_m\;\text{then }\eta_n^\otimes\ne\eta_m^\otimes.
\tag{$*_3$}\end{equation}
\begin{equation}
\eta_n^*\ne\eta_m^*\;\text{implies }\eta_n^\otimes\ne\eta_m^\otimes.
\tag{$*_4$}\end{equation}
Next, using again the Claim, choose $e_n^\otimes<\omega$ such that the
following condition holds true:
\begin{enumerate}
\item[(i)] $2\cdot e_n\leqslant e_n^\otimes$,
\item[(ii)] if $\eta_n^\otimes\ne \eta_m^\otimes$ then the sets
$\{\eta_n^\otimes\restriction l\colon e_n^\otimes\leqslant l<\omega\}$ and
$\{\eta_m^\otimes\restriction l\colon e_m^\otimes\leqslant l<\omega\}$ are
disjoint,
\item[(iii)] if $m<n<\omega$ and $m,n\in
Y^\otimes:=\{n<\omega\colon\eta_n^\otimes\notin\{\eta_l^\otimes\colon
l<n\}\}$, then $e_m^\otimes< e_n^\otimes$ and $Y\subseteq Y^\otimes$.
\end{enumerate}
Finally, for every $\eta\in{}^\omega 2$ we define:
$$\eta^{\text{even}}=\langle\eta(2i)\colon 2i<\text{lg}(\eta)\rangle,$$
$$u(\eta)=\{n\in Y^\otimes\colon(\exists
l<\omega)(e_n^\otimes<l<\omega\And \eta\restriction
l=\eta_n^\otimes\restriction l)\},$$
$$n_k(\eta)=\text{the k-th member of }u(\eta),$$
$$m_k(\eta)=\min\{m<\omega\colon \eta^{\text{even}}\restriction(m+1)\notin
C_{n_k}(\eta)\},$$
i.e. $m_k(\eta)$ is the smallest $m<\omega$ such that
$$\langle\eta(2i)\colon i<m+1\rangle\ne\eta^*_{n_k(\eta)}\restriction
(m+1).$$
Note that by the definition of $u(\eta)$ and $\eta^{\text{even}}$ and $m_k(\eta)$ we have
$$e_{n_k(\eta)}^\otimes<2\cdot m_k(\eta),$$
hence $m_k(\eta)>e_{n_k}(\eta)$. Clearly we also have:

\begin{enumerate}
\item[$(i)_1$] $u(\eta)$ is well defined if $\eta\in{}^{\omega}2$,
\item[$(ii)_2$] $n_k(\eta)$ is well defined if $k<|u(\eta)|$,
\item[$(iii)_3$] $m_k(\eta)$ is well defined if $k<|u(\eta)|$ and $\eta\ne
\eta^\otimes_{n_k}(\eta)$.
\end{enumerate}
Now we can define a function $\tau\colon
{}^{\omega}2\setminus\{\eta_n^\otimes\colon
n<\omega\}\to{}^{\omega\geqslant}2$ by the formula:
$$\tau(\eta)=f_{n_0(\eta)}(m_0(\eta))\:\sphat\:f_{n_1(\eta)}(m_1(\eta))\:\sphat\:\dots,$$
where, for $n<\omega$ $f_n$ is the function from the condition (3). 
\relax From the formula it easily follows that
$\tau(\eta)\in{}^{\omega\geqslant}2$ and it is well defined if
\mbox{$\eta\notin\{\eta_n^\otimes\colon n<\omega\}$} and moreover
$\tau(\eta)$ is infinite whenever $u(\eta)$ is infinite, as
$\langle\rangle\notin\text{Rang}(f_n)$.

To complete the proof of the theorem it remains to show that:

\begin{fact}[8]
$\Vdash_{\Bbb{P}}``\tau(\underline{\underline{r}}^\otimes)\;\text{is Cohen
over }V"$.  \end{fact}
To prove this fact we fix an open dense set $I\subseteq {}^{\omega>}2$ and
a $p\in\Bbb{P}$. Then $p\in A_n$ for some $n<\omega$. Let
$n^\otimes=\min\{m<\omega\colon\eta^\otimes_m=\eta^\otimes_n\}$. Clearly,
$n^\otimes\leqslant n$ and $n^\otimes\in Y^\otimes$. Then
$u(\eta^\otimes_n)$ is well defined and $n^\otimes\in u(\eta_n^\otimes)$;
in fact $n^\otimes$ is the last member of  $u(\eta_n^\otimes)$.
Let $k=|u(\eta_n^\otimes)|-1$, so $n_k(\eta_n^\otimes)=n^\otimes$. Also
$m_l(\eta_n^\otimes)$ is well defined and finite for $l<k$. Then we set
$$\nu=f_{n_0(\eta_n^\otimes)}(m_0(\eta^\otimes_n))\:\sphat\:
 f_{n_1(\eta_n^\otimes)}(m_1(\eta^\otimes_n))\:\sphat\dots\:\sphat\:
 f_{n_{k-1}(\eta_n^\otimes)}(m_{k-1}(\eta^\otimes_n)),$$
so if $k=0$ i.e.  $u(\eta_n^\otimes)$ the singleton then $\nu$ is the empty
sequence.

Clearly $\nu\in{}^{\omega>}2$. We can find  $p_1\in A_n$, $p\leqslant p_1$
such that
$$p_1\nVdash\neg``\underline{\underline{r}}^\otimes\restriction(e_n^\otimes+8)\unlhd\eta_n^\otimes\:".$$

Let $\varphi=``u(\eta_n^\otimes)$ is an initial segment of
$u(\underline{\underline{r}})^\otimes"$.  Clearly $$p_1\nVdash \neg``\varphi\;\text{and }l<k\;\text{implies
 }n_l(\underline{\underline{r}}^\otimes)=n_l(\eta_n^\otimes)\;\text{and
 }m_l(\underline{\underline{r}}^\otimes)=m_l(\eta_n^\otimes)"$$
 So however
 $$p_1\Vdash
 ``\underline{\underline{r^\otimes}}\ne\eta_n^\otimes\;\text{hence
 }\varphi\text{ implies } m_k(\underline{\underline{r}}^\otimes)\;\text{is well
defined}".$$
 Let
 $$Z=\{\rho\in{}^{\omega>}2\colon\neg(p_1\Vdash``f_{n_k}
 (\underline{\underline{r}}^\otimes)\ne\rho\text{ or }\neg\varphi")\};$$
equivalently $Z$ consists of all those $\rho\in{}^{\omega>}2$ for which
$p_1\nVdash\neg``\varphi\text{ and }f_{n_k}
 (\underline{\underline{r}}^\otimes)=\rho"$.
 So clearly it is enough to prove that $Z$ is a somewhere dense
 subset of ${}^{\omega>}2$. But
$$p_1\nVdash\neg``n_k(\underline{\underline{r}}^\otimes)=n^\otimes\text{
or }\neg\varphi",$$ so $Z=\{\rho\in{}^{\omega>}2\colon\neg(p_1\Vdash
``f_n^\otimes(m_k(\underline{\underline{r}}^\otimes))\ne\rho\text{ or
}\neg\varphi")\}$.  By the choice of $f_n^\otimes$ it is enough to prove
that:  $$B_0=\{m<\omega\colon \neg(p_1\Vdash``m_k(\underline{\underline{r}}
^\otimes)\ne m\text{ or }\neg\varphi")\}\in{\cal{D}}_n^\otimes.$$
But, by the condition ($*_4$), we have ${\cal{D}}_n^\otimes={\cal{D}}_n$.
So it is enough to prove that
$B_0\in{\cal{D}}_n$.
By the definition  of $m_k
(\underline{\underline{r}}^\otimes)$ this means, as
$n_k(\underline{\underline{r}}^\otimes)= n^\otimes$, that
$$\{m<\omega\colon \neg(p_1\Vdash``\underline{\underline{r}}
\restriction m\ne\eta_{n^\otimes}^*\restriction m\;\text{ or }
\underline{\underline{r}}\restriction(m+1)=\eta^*_{n^\otimes}
\restriction(m+1)\text{ or }\neg\varphi")\}\in{\cal{D}}_n,$$
but $\eta^*_n=\eta^*_{n^\otimes}$ as $\eta_n^*$ is defined from
$\eta_n^\otimes$. So this means that
$$\{m<\omega\colon(p_1\nVdash \neg(\underline{\underline{r}}
\restriction m=\eta_n^*\restriction m\;\text{and }
\underline{\underline{r}}\restriction (m+1)\ne\eta_n\restriction
(m+1) \text { and } \varphi)\}\in{\cal{D}}_n,$$
but $p_1\in A_n$, so this set belongs to ${\cal{D}}_n^0$ and
${\cal{D}}_n^0\subseteq{\cal{D}}_n$. So we are done. \end{pf}

Finally we prove that the theorem converse to Theorem 2 is also true, i.~e.
we shall show that whenever there exists a nowhere dense ultrafilter there
exists a $\sigma$-centered forcing $\Bbb{P}$ with the property that above
each element there are two incompatible ones and moreover $\Bbb{P}$ does
not add any Cohen real. To prove this fact we shall use some topological
methods.

Recall, a subalgebra $\Bbb{B}$ of a Boolean algebra $\Bbb{A}$ is
$\em{regular}$ whenever $\sup_{{\Bbb{A}}}X=1$ for every
$X\subseteq\Bbb{B}$ such that $\sup_{{\Bbb{B}}}X=1$.
The subalgebra $\Bbb{B}$ is regular iff the corresponding map of the Stone
spaces is semi-open, i.~e., the image of every non-empty clopen set  has
non-empty interior. Using Baumgartner's Principle that there exists
nowhere dense ultrafilters we construct a dense in itself separable
extremally disconnected compact space (= Stone space of a n
$\sigma$-complete Boolean algebra), which has no semi-open continuous maps
onto the Cantor set.

We use a topology on the set ${}^{\omega>}\omega=\bigcup\{{}^n\omega\colon
n<\omega\}$. If $s\in{}^{\omega>}\omega$ is a sequence of length $n$ and
$k\in\omega$, then $s\sphat\: k$ denotes the sequence of length $n+1$
extending $s$ in such a way that the $n$-th term is $k$. For a set
$A\subseteq \omega$ we set $s\sphat\: A=\{s\sphat\: k\colon k\in A\}$.
For a given ultrafilter $p\subseteq{\cal{P}}(\omega)$ we consider a
topology ${\cal{T}}_p$ on ${}^{\omega>}\omega$ given by the formula:
$$U\in{\cal{T}}_p\text{ iff for every }s\in U\text{ there exists }A\in
p\text{ such that }s\sphat\: A\subseteq U.$$
The set ${}^{\omega>}\omega$ equipped with the topology ${\cal{T}}_p$ we
denote $G_p$. The space $G_p$ is known to be Hausdorff and extremally
disconnected; see e.~g. Dow, Gubbi and Szymanski, (\cite {DGS}).
Hence the \v{C}ech-Stone extension $\beta G_p$ is an extremally
disconnected, compact, separable and dense in itself.

\begin{theor}
If there exists a nowhere dense ultrafilter then there exists a
$\sigma$-centered forcing $\Bbb{P}$ such that above every element of
$\Bbb{P}$ there are two incompatible ones and $\Bbb{P}$ does not add any
Cohen real.
\end{theor}

\begin{pf}
In virtue of a theorem of Silver it is enough to show that there exists a
$\sigma$-centered, complete, atomless Boolean algebra $\Bbb{B}$ such that
$\Bbb{B}$ does not contain any regular free subalgebra. For this goal we
shall use the topological space $G_p$ described above.
It remains to show that whenever $p$ is a nowhere dense ultrafilter and
$f\colon \beta G_p\to {}^\omega\{0,1\}$ is continuous, then there exists a
non-empty clopen set $H\subseteq\beta G_p$ such that $\inte
f(H)=\emptyset$.

First of all we notice that since $p$ is a nowhere dense ultrafilter, for
every $s\in{}^{\omega>}\omega$ there exists $A_s\in p$ such that
\begin{equation}\label{first}
\inte \cl f(s\sphat\:A_s)=\emptyset.
\end{equation}
In the sequel $L_n$ will denote the set of all sequences of length $n$,
i.~e., $L_n$ is the $n$-th level of the tree  ${}^{\omega>}\omega$. In
particular, $L_0=\{s_0\}$ is the empty sequence. By induction we
define a sequence of sets $\{U_n\colon n<\omega\}$ such that $U_n\subseteq
L_n$ for every $n<\omega$ and, moreover
\begin{equation}\label{second}
\inte \cl f(U_n)=\emptyset,
\end{equation}
\begin{equation}\label{third}
\text{for every }s\in U_n\text{ there exists }A\in p\text{ such that
}s\sphat\: A\subseteq U_{n+1}.
\end{equation}
We set $U_0=\{s_0\}$ and $U_1=s_0\:\sphat\: A_{s_0}$. Assume $U_n$ is
defined, say $U_n=\{s_k\colon k<\omega\}$. Then by continuity of $f$ and
the condition (\ref{first}) we can choose $A_k\in p$ in such a way that
$\inte \cl f(s_k\:\sphat\: A_k)=\emptyset$ and moreover, the diameter of
$\cl f(s_k\:\sphat\: A_k)$ is not greater than~$\frac{1}{k}$. Clearly,
$s_k$ is an accumulation point of $s_k\:\sphat\: A_k$, because $A_k\in p$.
Hence, for every $k<\omega$ we get
$$\cl f(s_k\:\sphat\:A_k)\cap \cl f(U_n)\ne\emptyset.$$
Therefore, since diameters of the sets $\cl f(s_k\:\sphat\:A_k)$ tends
to zero, the set of accumulation points of the set $\bigcup\{cl
f(s_k\:\sphat\: A_k)\colon k<\omega\}$ is contained in $\cl f(U_n)$. 
Indeed, every $\varepsilon$-neighbourhood of the set $\cl f(U_n)$ has to
contain all but infinitely many sets of the form $\cl f(s_k\sphat A_k)$. 
Now
we set
$$U_{n+1}=\bigcup\{s_k\:\sphat\:A_k\colon k<\omega\}$$
and observe that
$$\cl f(U_{n+1})\subseteq \cl f(U_n)\cup\bigcup\{\cl
f(s_k\:\sphat\:A_k)\colon k<\omega\}.$$
So the set 
$\cl f(U_n)\cup\bigcup\{\cl f(s_k\:\sphat\:A_k)\colon k<\omega\}$ is
closed.
Thus the set $\cl f(U_{n+1})$ is nowhere dense, because both the sets $\cl
f(U_n)$ and $\cl f(s_k\:\sphat\:A_k)$ for $k<\omega$, are nowhere dense
and the set on the right hand side is closed.

By the condition (\ref{second}) and the Baire Category Theorem, there
 exists a dense set
$$\{x_n\colon n<\omega\}\subseteq {}^\omega\{0,1\}\setminus\bigcup\{\cl
f(U_n)\colon n<\omega\}.$$
In particular, for every $n,k<\omega$ we have
$$f^{-1}(\{x_n\})\cap \cl U_k=\emptyset,$$
where ``$\cl$" denotes here the closure in $\beta G_p$. Now, for every
 $n<\omega$ we choose a clopen set $V_n\subseteq \beta G_p$ such that
\begin{equation}\label{fourth}
f^{-1}(\{x_n\})\subseteq V_n\subseteq \beta G_p\setminus(\cl
 U_0\cup\dots\cup U_n).
 \end{equation}
By induction we construct a sequence $\{W_n\colon n<\omega\}$ such
that the following conditions hold true: 
\begin{equation}\label{fifth}
W_n\subseteq U_n\text{ for }n<\omega\text{ and }W_0=U_0
\end{equation}

for every $ s\in W_n $ there exists $ B_s\in p $ such that 

\begin{equation}\label{sixth}
s\:\sphat\:B_s\subseteq U\setminus(V_0\cup\dots\cup V_n),
\end{equation}

\begin{equation}\label{seventh}
W_{n+1}=\bigcup\{s\:\sphat\:B_s\colon s\in W_n\}.
\end{equation}

Assume the sets $W_0,\dots,W_n$ are defined in such a way that
(\ref{fifth}), (\ref{sixth}) and (\ref{seventh}) are satisfied.  Then we
have in particular 
$$W_n\subseteq U_n\setminus(V_0\cup\dots\cup V_{n-1});$$
by the condition (\ref{fourth}) we also have 
$$U_n\subseteq \beta G_p\setminus V_n.$$
Hence we get $W_n\subseteq U_n\setminus(V_0\cup\dots\cup V_n)$. Since the
set $U_n\setminus (V_0\cup\dots\cup V_n)$ is open, for every $s\in W_n$ we
can choose $B_s\in p$ such that $s\:\sphat\:B_s\subseteq
U_n\setminus(V_0\cup\dots\cup V_n)$. Then it is enough to set
$W_{n+1}=\bigcup\{s\:\sphat\: B_s\colon s\in W_n\}$.

Clearly the set $W=\bigcup\{W_n\colon n<\omega\}$ is open in $G_p$ and
$W\cap V_n=\emptyset$ for every $n<\omega$. Indeed, if $m>n$, then
$W_m\cap V_n=\emptyset$ by the conditions (\ref{sixth}) and
(\ref{seventh}), whereas for $m\leqslant n$, $W_m\cap V_n=\emptyset$
because $W_m\subseteq U_m$ and $U_m\cap V_n=\emptyset$ by the condition
(\ref{fourth}). Since $V_n$ is a clopen in $\beta G_p$ we also have
$$\cl W\cap V_n=\emptyset$$
for every $n<\omega$. Since $\beta G_p$ is extremally disconnected, $\cl
W$ is clopen subset of $\beta G_p$ and, by the last equality and condition
(\ref{fourth}) we get
$$f(\cl W)\cap\{x_n\colon n<\omega\}=\emptyset.$$
Therefore $f(\cl W)$ is nowhere dense, because $\{x_n\colon n<\omega\}$ is
dense in ${}^\omega\{0,1\}$, which completes the proof.
\end{pf}

\end{document}